\documentclass[10pt]{article}

\usepackage[outerbars]{changebar} 
\usepackage{color} 
\usepackage{comment} 
\usepackage{amssymb}
\usepackage{amsmath}
\usepackage{ulem}
\usepackage{tikz-cd}

\def\AnnENS{Ann.~\'Ec.~Norm.~}

\def \ccomma{\raise 2pt\hbox{,}} 
\def \D {\hbox{d}}

\def\C{\mathbb C}
\def\S{\mathcal S}

\def \NDLR#1 {\hfill\break\noindent NDLR \textbf{#1}\hfill\break\noindent}

\begin{document}

\title{Equivalence of two canonical lists of elliptic functions}

\author{Robert Conte
{}\\
\\ 1. Universit\'e Paris-Saclay, ENS Paris-Saclay, CNRS,
\\    Centre Borelli, F-91190 Gif-sur-Yvette, France
\\
\\ 2. Department of mathematics, The University of Hong Kong,
\\ Pokfulam road, Hong Kong
\\    E-mail Robert.Conte@cea.fr, ORCID https://orcid.org/0000-0002-1840-5095
{}\\
 Liangwen Liao
{}\\
Department of Mathematics, Nanjing University, Nanjing, China
\\    E-mail maliao@nju.edu.cn,   ORCID https://orcid.org/0000-0002-3337-5146
{}\\
 Chengfa Wu
{}\\
Institute for Advanced Study, Shenzhen University, Shenzhen, China
\\    E-mail cfwu@szu.edu.cn, ORCID https://orcid.org/0000-0003-1697-4654
{}\\
}

\date{}

\maketitle

\hfill 

{\vglue -10.0 truemm}
{\vskip -10.0 truemm}


\begin{abstract}
We establish a one-to-one correspondence between,
on one hand
the four types of trans\-cendental meromorphic solutions
of the autonomous Schwarzian differential
equations which are elliptic,
on the other hand
the four binomial equations of Briot and Bouquet possessing an elliptic solution.
\end{abstract}




\baselineskip=12truept 

\vfill\eject



When looking for first order ordinary differential equations (ODEs)
whose general solution is singlevalued,
Briot and Bouquet   \cite{BriotBouquet,Hille}  enumerated and integrated all
first order binomial equations having such a property,
\begin{eqnarray} \label{first-order ODE}
& & \left(\frac{\D u}{\D z}\right)^k = R(u),
\label{eqBB-Binomial}
\end{eqnarray}
with $k$ a positive integer and $R$ a rational function.
\textit{Modulo} homographic transformations (which leave (\ref{eqBB-Binomial}) form-invariant),
only four of them are elliptic,
and their general solution is a homographic transform of, respectively,
$\wp$, $\wp'$, $\wp^2$, $\wp^3$,
where $\wp$ denotes the elliptic function of Weierstrass.


{}From another point of view,
namely the growth at infinity of meromorphic solutions of ODEs,
another class of ODEs also form-invariant under homographic transformations
has been thoroughly studied,
this class is made of the autonomous Schwarzian differential equations
\begin{eqnarray}
& &
\left\lbrace u;z \right\rbrace^p=R(u),
\left\lbrace u;z \right\rbrace=\frac{u'''}{u'}-\frac{3}{2} \left(\frac{u''}{u'}\right)^2,
\label{eqOrder3-degree1-Schwarzian-autonomous}
\end{eqnarray}
with $p$ a positive integer and $R$ a rational function.

Ishizaki \cite{Ishizaki-1991} enumerated all ODEs (\ref{eqOrder3-degree1-Schwarzian-autonomous})
possessing particular transcendental meromorphic solutions.
\textit{Modulo} homographies, there are only six of them,
whose canonical form with an equal number
(counting multiplicity) of poles and zeroes in $R$ is,
\begin{eqnarray}
& &
\left\lbrace
\begin{array}{ll}
\displaystyle{
\left\lbrace u;z \right\rbrace  =c
\frac{(u-\sigma_1)  (u-\sigma_2)(u-\sigma_3)(u-\sigma_4)}{(u-\tau_1)  (u-\tau_2)(u-\tau_3)(u-\tau_4)}\ccomma  
}\\ \displaystyle{
\left\lbrace u;z \right\rbrace^3=c
\frac{(u-\sigma_1)^3(u-\sigma_2)^3}                      {(u-\tau_1)^2(u-\tau_2)^2(u-\tau_3)^2}\ccomma       
}\\ \displaystyle{
\left\lbrace u;z \right\rbrace^2=c
\frac{(u-\sigma_1)^2(u-\sigma_2)^2}                      {(u-\tau_1)^2(u-\tau_2)(u-\tau_3)}\ccomma          
}\\ \displaystyle{
\left\lbrace u;z \right\rbrace^3=c
\frac{(u-\sigma_1)^3(u-\sigma_2)^3}                      {(u-\tau_1)^3(u-\tau_2)^2(u-\tau_3)}\ccomma        
}\\ \displaystyle{
\left\lbrace u;z \right\rbrace  =c
\frac{(u-\sigma_1)  (u-\sigma_2)}                        {(u-\tau_1)  (u-\tau_2)}\ccomma                
}\\ \displaystyle{
\left\lbrace u;z \right\rbrace  =c,                                                                                                     
}
\end{array}
\right.
\label{eqList6}
\end{eqnarray}
where $c$, $\tau_j$, $\sigma_j$, $j=1,\dots,4$ are complex constants,
the $\tau_j$'s being distinct and
the $\sigma_j$'s depending on the $\tau_j$'s.

Then two of us \cite{LiaoWu-Schwarzian-ODEs} found explicitly
all the transcendental meromorphic solutions of Eqs.~(\ref{eqList6}),
except those solutions of Eq.~(\ref{eqList6})$_{5}$ that have no Picard exceptional values.
One of the key points in the derivation of these solutions is to show that each of them
obeys a first order differential equation of the form \eqref{first-order ODE}
by applying Nevanlinna theory \cite{LaineBook} on the value distribution of meromorphic functions
and the Laurent expansions around poles.
Nevertheless, for Eq.~(\ref{eqList6})$_{5}$, the argument in \cite{LiaoWu-Schwarzian-ODEs}
only works when $u$ possesses at least one Picard exceptional value.
Then transcendental meromorphic solutions of Eqs.~(\ref{eqList6}) can be obtained
by solving the corresponding first order differential equations.

We remark that all transcendental meromorphic solutions of Eqs.~(\ref{eqList6})  given in \cite{LiaoWu-Schwarzian-ODEs}  belong to class {\it $W$},
which consists of elliptic functions and their degeneracies,
i.e.~functions $R(z)$ or $R(e^{az}), a \in \C$, where $R$ is a rational function.
As shown by Eremenko \cite{EremenkoKS}, functions in class {\it $W$}
constitute all transcendental meromorphic solutions of a class of
autonomous algebraic differential equations \cite{ConteNgWong-RSH,EremenkoKS,NgWu2019-Loewy}.
It turns out that Ermenko's method can be further generalized to a wider class
of differential equations \cite{CNW-boundary-layer,Demina-Two-top-degree}.

For the last equation of (\ref{eqList6}), all the transcendental meromorphic solutions
are the three-parameter homographic transforms of $e^{a z}$
(or equivalently $\coth(a z/2)$), $c=-a^2/2$.

For the first four equations in (\ref{eqList6}),
all the transcendental meromorphic solutions are a one-parameter elliptic function.
Since the solutions  $u(z)$ of Theorems 5 and 6 in Ref.~\cite{LiaoWu-Schwarzian-ODEs}
obey the respective ODEs
\begin{eqnarray}
& & \hbox{(Thm 5) } \left(\frac{\D u}{\D z}\right)^3 = L^3 u^2 (u-1)^2 (u+1)^2,
\\
& & \hbox{(Thm 6) } \left(\frac{\D u}{\D z}\right)^4 = L^4 u^2 (u-1)^3 (u+1)^3,
\end{eqnarray}
they admit the general solutions
\begin{eqnarray}
& & \hbox{(Thm 5) } u(z)= - \frac{2 L^3}{27 \wp'(z-z_0)}\ccomma
\\
& & \hbox{(Thm 6) } u(z)= 1 + \frac{L^4}{32 (\wp^2(z-z_0) - L^4/64)}\ccomma
\end{eqnarray}
with $z_0$ arbitrary.

Therefore all the transcendental meromorphic solutions of the first four equations in the list (\ref{eqList6})
are the most general homographic transforms of, respectively,
$\wp(z-z_0)$, $\wp'(z-z_0)$, $\wp^2(z-z_0)$, $\wp^3(z-z_0)$.
This establishes a one-to-one correspondence with the above recalled result of Briot and Bouquet.

Choosing the identity for the above mentioned homographic transformations,
one deduces a second canonical form of equations (\ref{eqList6}),
\begin{equation}
\left\{ \begin{array}{lll}
\left\lbrace u;z \right\rbrace   & = \displaystyle{- \frac{3}{8} \frac{(4 u^2+g_2)^2 +32 g_3 u}{4 u^3-g_2 u -g_3}} & u=\wp(z,g_2,g_3) \\
\left\lbrace u;z \right\rbrace^3 & =\displaystyle{-16 \frac{(u^2-3 g_3)^3}{(u^2+g_3)^2}}                           & u=\wp'(z,0,g_3) \\
\left\lbrace u;z \right\rbrace^2 & =\displaystyle{\frac{9 (80 u^2-24 g_2 u + 5 g_2^2)}{64 u (4 u - g_2)^2}}        & u=\wp^2(z,g_2,0) \\
\left\lbrace u;z \right\rbrace^3 & =\displaystyle{-8 \frac{(35 u^2 - 10 g_3 u + 2 g_3^2)^3}{u^2 (4 u - g_3)^3}} & u=\wp^3(z,0,g_3) \\
\left\lbrace u;z \right\rbrace & =\displaystyle{\frac{a^2 (2 u^2-a^2)}{u^2+a^2}} & u=\displaystyle{\frac{a}{\sinh(a z)}} \\
\left\lbrace u;z \right\rbrace  &  =\displaystyle{-\frac{a^2}{2}} & u=\displaystyle{\frac{c_1 e^{a z}+c_2}{c_3 e^{a z}+c_4}}\cdot
\end{array}\right.
\label{eqList6bis}
\end{equation}

\textit{Remark}.
The first five equations (\ref{eqList6bis}) have no other singlevalued solutions
than the above mentioned one-parameter meromorphic functions.
Indeed,
each of them admits at least one Laurent series
whose Fuchs indices \cite{CMBook} are $j=-1$ and two irrational numbers,
respectively
\begin{equation}
%
\left\{ \begin{array}{lll}
u  & =    z^{-2} + \dots & \hbox{indices } (j+1)(j^2-j-3)=0, \\ 
u  & = -2 z^{-3} + \dots & \hbox{indices } (j+1)(j^2-j-8)=0, \\ 
u  & =    z^{-4} + \dots & \hbox{indices } (j+1)(j^2-j-15)=0, \\ 
u  & =    z^{-6} + \dots & \hbox{indices } (j+1)(j^2-j-35)=0, \\ 
\displaystyle{\frac{1}{u-i a}} & = \displaystyle{\frac{2 i}{a^3}} z^{-2} + \dots & \hbox{indices } (j+1)(j^2-j-3)=0, \\ 
\end{array}\right.
\label{eqList5}
\end{equation}
thus lowering by two the maximal number (three) of arbitrary constants
in any singlevalued solution.


\vfill\eject


\end{document}